\documentclass [12pt]{article}
\usepackage{graphicx,amssymb,amsfonts,latexsym,amsmath,amsthm,times}
\usepackage{epsfig}
\usepackage{fancyhdr}
\usepackage{color}
\usepackage[english]{babel}
\numberwithin{equation}{section}

\begin{document}

%
%

\thispagestyle{plain}

\vspace*{2cm} \normalsize \centerline{\Large \bf Convergence of
Time-Dependent Turing Structures}

\medskip

\centerline{\Large \bf to a Stationary Solution}

\vspace*{1cm}

\centerline{A. G. Ramm$^1$, V. Volpert$^2$ \footnote{Corresponding
author. E-mail: volpert@math.univ-lyon1.fr}}

\vspace*{0.5cm}


\centerline{$^1$ Mathematics Department, Kansas State University,
Manhattan, KS 66506-2602, USA}
\centerline{email: ramm@math.ksu.edu}

 \centerline{$^2$ Institut Camille Jordan, UMR 5208 CNRS, University Lyon
1} \centerline{69622 Villeurbanne, France}


\vspace*{1cm}

\noindent {\bf Abstract.} Stability of stationary solutions of
parabolic equations is conventionally studied by linear stability
analysis, Lyapunov functions or lower and upper functions. We
discuss here another approach based on differential inequalities
written for the $L^2$ norm of the solution. This method is
appropriate for the equations with time dependent coefficients.
It yields new results and is applicable when the usual linearization
method is not applicable.

\vspace*{0.5cm}




\vspace*{1cm}


 \section{Formulation of the problem}
Large-time behavior of solutions to differential equations has
been discussed in many publications, see, for example, \cite{dk},
\cite{Kl}, \cite{Sho}, \cite{Zheng}. First, one has to establish
the global existence of the solution. This is done in most cases
by establishing an a priori estimate which implies boundedness of
the solutions for all times. The usual approach to Lyapunov
stability of solutions is to linearize the problem and prove that
the spectrum of the linearized operator lies strictly in the left
half-plane of the complex plane.

In recent papers \cite{R4}, \cite{R1}, \cite{R2}, \cite{R3},  a
novel approach to the stability and long-time behavior of
solutions to abstract differential equations is developed. This
approach is applied here to nonlinear systems of interest in
biology.

 Consider the semilinear parabolic system of equations

 \begin{equation}
 \label{m1}
 \frac{\partial u}{\partial t} = D(t) \Delta u  + F(u,x,t)
 \end{equation}
 in a bounded domain $\Omega \subset \mathbb R^{\mathcal{M}}$ with a
 sufficiently smooth boundary and with the homogeneous Dirichlet or Neumann
 boundary condition

 \begin{equation}
 \label{m2}
 u |_{\partial \Omega} = 0 \;\;\; {\rm or} \;\;\;
 \frac{\partial u}{\partial n}  |_{\partial \Omega} = 0 ,
 \end{equation}
 and the initial condition

 \begin{equation}
 \label{m3}
 u(x,0) = u_0(x) .
 \end{equation}
 Here $u=(u_1,...,u_n), F=(F_1,...,F_n)$, $D$ is a diagonal matrix with
positive diagonal elements $d_i=d_i(t)$,  which can depend on $t$,
and
 \begin{equation}
 \label{m4}
 F(0,x,t) = 0 , \;\; \forall x \in \Omega, t \geq 0 .
 \end{equation}
 The vector-function $F$ is assumed to satisfy the estimates
 \begin{equation}
 \label{m5}
\sup_{u, x\in R^{\mathcal{M}}, t\geq 0} |F(u,x,t)|\le M_1,
\end{equation}
\begin{equation}
 \label{m6}
|F(v,x,t)-F(w,x,s)|\le c_F(|t-s|+|v-w|),
\end{equation}
where $c_F>0$ is a constant independent of $v,w,x,t$, and $F$ is
continuous with respect to $x$.
 Under these conditions, $u=0$ is a stationary solution of problem
 (\ref{m1}), (\ref{m2}). In the examples considered below we assume that
$\mathcal{M}\le 3$.

Consider the operator linearized about
 this solution:

 $$ L_t v = D \Delta v + F_u'(0,x,t) v $$
 acting in the Hilbert space $L^2(\Omega)$ with the domain

 $$ {\mathcal D} = \{ u \in H^2(\Omega) , \; u |_{\partial \Omega} = 0 \;\;
 {\rm or} \;\;
 \frac{\partial u}{\partial n} |_{\partial \Omega} = 0 \} . $$
 Here $t$ is considered as a parameter.

 Suppose that the spectrum of the operator $L_t$ is located in the
 half-plane Re $\lambda \leq \sigma(t)$. If
 $$ \sigma(t) \leq \sigma_0 < 0 , \;\;\; t \geq 0 , $$
 then solution of problem (\ref{m1})-(\ref{m3}) converges to the
 stationary solution $u=0$.  This means that the stationary
(equilibrium)
solution is exponentially stable, i.e.,  the solutions with sufficiently
small initial data converge to the stationary solution $u=0$ at an
exponential rate.

A proof of this assertion is well  known in the case of the abstract
evolution problem of the type
\begin{equation}
 \label{evo}
\dot{u}=Au +B(t,u),\quad u(0)=u_0,
\end{equation}
where $A$ is a linear bounded operator in a Banach space,
with the spectrum that lies in the half-plane Re$z\leq \sigma_0<0$, and
$B(t,u)$ is a nonlinear operator satisfying the assumption
$$||B(t,u)||\leq c_0(t)||u||^p, \quad p>1,$$
where $c_0(t)$ satisfies
a suitable smallness assumption (see, e.g., \cite{dk}).

In \cite{dk}, Theorem I.4.1, the following result is proved: if $A$
is a bounded operator in a Banach space, then there exists the limit
$$\kappa:=\lim_{t\to \infty}\frac{\ln \|e^{At}\|}{t}=\max\{Re
\lambda|\lambda\in
\sigma(A)\},$$
where $\sigma(A)$ is the spectrum of $A$. Therefore, if all the
 solutions to the Cauchy problem
$$\dot{u}=Au,\quad u(0)=u_0,   $$
decay exponentially fast to zero, then the spectrum of $A$ lies
in the half-plane Re$z\leq -\kappa$, $\kappa>0$, and vice versa.
One should have in mind that if $A$ is a bounded linear operator in a
Hilbert space $H$,  such that Re$A\le -\kappa$, $\kappa>0$,
i.e., Re$(Au,u)\le -\kappa\|u\|^2$, then the spectrum of $A$ lies
in the half-plane  Re$z\leq -\kappa$, but the converse of this statement
is false if $\dim H>1$: even in two-dimensional Hilbert space one can
give an example of $A$ with the spectrum lying in the half-plane
Re$z\leq -\kappa<0$, for which the inequality Re$(Au,u)\le -\kappa\|u\|^2$
does not hold.    For instance, let
 $$ A = \left( \begin{array}{ccc}
 -1 &  & 3 \\
 0 &  & -1 \end{array} \right) , $$

The spectrum of this $A$ consists of negative eigenvalue $\lambda=-1$.
The quadratic form for real-valued $u_1$ and $u_2$ is
Re$(Au,u)=-u_1^2-u_2^2+3u_1u_2$. This quadratic form is {\it not}
negative-definite.

 If the spectrum of $A$ does not lie strictly in the
left half-plane of the complex plane, or $\sigma(t) \to 0$ as $t \to
 \infty$, then the assertion about exponential rate of convergence to zero
of the solutions to the Cauchy problem (\ref{evo}) is not valid, in
general, and the
Lyapunov stability problem cannot be solved by a study of the linearized
problem.

 In this work we  study this, more difficult, case, and use
a new technical tool for such a study, see Lemma 2.1.
 Let us emphasize that
 $\sigma (t)$ will not necessarily be assumed negative in this paper
(see \cite{R3}).


 \section{Convergence of solutions}

 In what follows we assume that
 $F(u,x,t)$ satisfies assumptions made in Section 1, see (\ref{m5})
and (\ref{m6}).
 The initial and the boundary conditions satisfy the compatibility
 conditions,
 $u_0(x)=0$ on $\partial \Omega$ for the Dirichlet and
 $\frac{\partial u}{\partial n}=0$ on $\partial \Omega$
 for the Neumann boundary condition.
 Under these (and some additional) conditions (see, e.g.,  \cite{LSU})
there exists a classical solution of problem
 (\ref{m1})-(\ref{m3}).

 Let us assume that
 $$ F(u,x,t) = A(x,t) u + B(u,x,t) , $$
 where

 \begin{equation}
 \label{c1}
 {\rm Re} \; (A(x,t) u, u) \leq - \gamma(t) |u|^2 ,  \;\;\; \forall x \in
 \Omega, \; t \geq 0,
 \end{equation}
 and

 \begin{equation}
 \label{c2}
 |B(u,x,t)| \leq c_0(t) |u|^p ,  \;\;\; \forall x \in
 \Omega, \; t \geq 0, \quad p>1.
 \end{equation}
 Here $( , )$ denotes the inner product in $\mathbb R^3$, and
$|u|^2=\sum_{j=1}^n |u_i|^2$.
 We  assume that the diagonal elements $d_i(t)$ of the matrix $D(t)$ of
the diffusion
 coefficients satisfy the estimates

 \begin{equation}
 \label{c3}
 d_i(t) \geq d(t) , \;\;\; i=1,...,n, \;\; t \geq 0 ,
 \end{equation}
 where $d(t)$ is a positive function. The assumptions about the behavior
of $d(t)$ for large $t$ will be specified below, in the formulation of
Theorems 3.1-3.3.

 \medskip

 Let $g(t):= \|u(\cdot,t)\|$, where $\| \cdot \|$ denotes the
$L^2(\Omega)$ norm. We will also use the space $L^\infty(\Omega)$
with the norm $\| \cdot \|_{\infty}$, and the usual Sobolev space
$H^2(\Omega)$ with the norm $\| \cdot \|_{H^2(\Omega)}$.

Multiplying equation (\ref{m1}) by $u$ and integrating, we
 obtain, taking into account (\ref{c1})-(\ref{c3}):

 \begin{equation}
 \label{c5}
  g \dot g \leq - d(t) \| \nabla u\|^2 - \gamma(t) g^2 + c_0(t)
  \int_\Omega |u|^{p+1} dx .
 \end{equation}
 In the case of the Dirichlet boundary condition, we use the
 Poincar\'e inequality

 \begin{equation}
 \label{c5a}
 c(\Omega) \int_\Omega |u|^2 dx \leq \int_\Omega |\nabla u|^2 dx ,
 \end{equation}
 where $c(\Omega)$ is a positive constant which depends on the domain.
The optimal (maximal possible) value of $c(\Omega)$ is equal
to the first eigenvalue $\lambda_1$ of the Dirichlet Laplacian in
$\Omega$.

In the case of the Neumann
 boundary condition, we put $c(\Omega)=0$.

 Using the following  multiplicative inequality (see, e.g., \cite{Bur},
p.193):
 $$ \|u\|_\infty \le c
 \|u\|_{L^2(\Omega)}^{1/4} \|u\|_{H^2(\Omega)}^{3/4} , $$
where the constant $c>0$ is independent of $u$,
 we obtain

 $$ \int_\Omega |u|^{p+1} dx \le g^2 ||u||_\infty^{p-1} \leq
  c^{p-1} \|u\|_{H^2(\Omega)}^{3(p-1)/4} g^{(p+7)/4}. $$
 From this inequality, (\ref{c5}) and (\ref{c5a}) we obtain
 \begin{equation}
 \label{c6}
  \dot g \leq - (d(t) c(\Omega) + \gamma(t)) g + c_0(t)
  c^{p-1} \|u\|_{H^2(\Omega)}^{3(p-1)/4} g^{(p+3)/4} .
 \end{equation}

 \noindent
  It is known that under our assumptions the $H^2$ norm of the
solution is bounded (see \cite{f}, Theorem 16.1, p.170,  and  Section
4).

Define
 \begin{equation}
 \label{alpha}
 \sigma(t) := d(t) c(\Omega) + \gamma(t), \;\;\;
 \alpha(t) := c_0(t) c^{p-1} \|u\|_{H^2(\Omega)}^{3(p-1)/4}, \;\;\;
 q := \frac{p+3}{4} .
 \end{equation}
 Then (\ref{c6}) can be written as

 \begin{equation}
 \label{c7}
  \dot g \leq - \sigma(t) g + \alpha(t) g^q , \qquad g(t)\ge 0,
 \end{equation}
where $q>1$ because $p>1$.

 Assume that $\sigma(t)$ and $\alpha(t)\ge 0$ are continuous functions
defined on $[0,\infty)$.

 We will use in Section 3 the following basic result from \cite{R1},
where more general results are obtained (see
also \cite{R3}):

 {\bf Lemma 2.1}. {\it If there exists
 a function $\mu(t)>0$, defined on $[0,\infty)$,  such that

 \begin{equation}
 \label{c8}
  \alpha(t)  \leq  \mu^{q-1}(t) \left( \sigma(t) - \frac{\dot
\mu(t)}{\mu(t)}
  \right) , \;\; t \geq 0
 \end{equation}
 and

 \begin{equation}
 \label{c9}
   \mu(0) g(0) \leq 1 ,
 \end{equation}
 then $g(t)$ exists for all $t\ge 0$, and

 \begin{equation}
 \label{c10}
  0 \leq g(t) \leq \frac{1}{\mu(t)} , \;\; \forall t \geq 0 .
 \end{equation}
}

Note that if $\lim_{t\to \infty}\mu(t)=\infty$, then
estimate (\ref{c10}) implies that $\lim_{t\to \infty}g(t)=0$.
The function $\sigma(t) $ in lemma 2.1 is not necessarily positive.


 \section{Applications}

A relatively general class of abstract differential equations
for which our method is applicable is described by the equations
of the form
$$\dot{u}=A(t)u+G(t,u)+f(t), \quad u(0)=u_0,$$
where $A$ is a linear operator in a Hilbert space $H$, $G$ is a
nonlinear operator in  $H$, and $f$ is a given function with
values in  $H$. The following assumptions allow one to use our
approach: Re$(A(t)u,u)\le -\gamma(t)||u||^2$, $||G(t,u)||\le
a(t,g)$, $g:=||u(t)||$, $||f(t)||\le \beta(t)$, where the
functions $\gamma$,$a(t,g)$ and $\beta$ satisfy some assumptions
that are detailed in \cite{R4}, \cite{R1}.

 In this section we will apply the results obtained above to
 reaction-diffusion system (\ref{m1}) with time dependent
 coefficients. In particular, in the case where the diffusion
 coefficients converge to zero and conventional results on
 stability of stationary solutions are not applicable.

 \newpage

 \subsection{Convergence with various rates}

 \bigskip

 {\bf Exponential rate of convergence.}

 In order to make clear our method for a study of the large-time behavior
of the solution to problem
 (\ref{m1})-(\ref{m3}), let us consider first  a single equation and the
 Dirichlet boundary condition.

 Specifically,  consider the following example:

 $$ A(x,t) \equiv a_0 > 0, \;\;\; D(t) \equiv d_0 > 0 , \;\;\;
 c_0(t) \equiv c_0, $$
 where $a_0, d_0$ and $c_0$ are some constants, and $p=2$ in (\ref{c2}).
Then $\gamma(t)=-a_0$, and
 $$ \sigma = d_0 c(\Omega) - a_0 , \;\;\; q=\frac54 \; . $$
 If $\sigma > 0$, that is, if

 \begin{equation}
 \label{ap1}
 \frac{a_0}{d_0} < c(\Omega) ,
 \end{equation}
 then we choose

 $$ \mu(t) = \mu_0 e^{\nu t} , $$
 where $\mu_0$ and $ \nu$ are positive constants.

Let us formulate sufficient conditions for assumptions (\ref{c8}) and
(\ref{c9}) to be satisfied. If these assumptions are  satisfied, then
inequality (\ref{c10}) yields an exponential rate of decay of the function
$g(t)$, and, therefore, of the solution $u(t)$ to zero.

 This assertion can be explained in terms of the exponential stability
in the sense of Lyapunov  of
the  solution $u=0$ to the problem  (\ref{m1})-(\ref{m3}). Namely,
consider the
 problem, linearized about the zero solution.
 The principal eigenvalue of the linearized problem

 $$ d_0 \Delta u + a_0 u = \lambda u , \;\;\; u |_{\partial
 \Omega} = 0 $$
 becomes negative if

 \begin{equation}
 \label{ap2}
 \frac{a_0}{d_0} < c(\Omega),
 \end{equation}
 where the constant $c(\Omega)$ is from the Poincar\'e inequality
(\ref{c5a}).
 This is  condition (\ref{ap1}).

To satisfy assumption (\ref{c9}) one may choose
$$\mu_0=g(0)^{-1}.$$
One may assume that $g(0)\neq 0$, because otherwise the solution
is zero by the uniqueness theorem that holds under our assumptions.

To satisfy assumption (\ref{c8}) it is sufficient to assume that
\begin{equation}
 \label{m7}
 C c_0(t)\leq g(0)^{-\frac 1 4}e^{{\frac 1 4}\nu t}[\sigma(t)-\nu],
\end{equation}
where we took into account that $q-1=\frac 1 4$ and denoted by $C$
the constant $c^{p-1}||u||^{3/4}_{H^2(\Omega)}$. In Section 4 it is proved
that the norm $||u||_{H^2(\Omega)}$ can be estimated from above by a
constant independent of $u$.
If one chooses $\nu=0.5 \sigma =(d_0\lambda_1 -a_0)/2$, then
inequality (\ref{m7}) holds provided that
\begin{equation}
 \label{m8}
c_0(t)\le 0.5 g(0)^{-\frac 1 4}C^{-1}e^{{\frac 1 8}\sigma t}\sigma.
\end{equation}
Thus, condition  (\ref{m8}) is sufficient for the assumption (\ref{c8})
to be satisfied. Condition (\ref{m8}) holds for any fixed $g(0)$ if
$c_0(t)$ is sufficiently small. It holds for a fixed $c_0(t)$
if $g(0)$ is sufficiently small.

 We have proved the following result.

 \medskip

 \noindent
 {\bf Theorem 3.1.}
 {\it Let the function $\sigma(t)$, defined in (\ref{alpha}), satisfy
 the inequality
$$ \sigma(t) \geq \sigma_0 >0,\quad  \forall t\ge 0, $$
 where $\sigma_0$
 is a constant. Choose $\mu(t) = \mu_0 \exp(\nu t)$, $\mu_0=g(0)^{-1}$,
$\nu =0.5 \sigma_0$. If condition  (\ref{m8}) holds,
 then the $L^2$ norm of the solution $u(x,t)$ to problem
 (\ref{m1})-(\ref{m3}) with the Dirichlet boundary condition
 satisfies the estimate
\begin{equation}
\label{in2}
  \|u(\cdot,t)\| \; \leq \; g(0) e^{-0.5 \sigma_0 t} \; , \;\;\; \forall t
 \geq 0 .
\end{equation}
}
 \medskip

 \noindent
 The conclusion of this theorem follows from Lemma 2.1, see estimate
(\ref{c10}).

The method for  estimating the large time behavior of solutions to
evolution problems, that was used in the proof of Theorem 3.1 is
easy to apply in many problems.

 The assumptions  of Theorem 3.1 do not explicitly
 require that the spectrum of the linearized problem lies
 in the open left half-plane of the complex plane. However, the
exponential rate of decay
 of the solution suggests that this is the case (see \cite{dk},
p.42, p.51).


  \bigskip

 {\bf Convergence at a power rate.}

Let us consider
 problem (\ref{m1})-(\ref{m3}) with the Dirichlet  boundary condition.
 Let us assume that

 \begin{equation}
 \label{ap3a}
 d(t) = \frac{d_0}{t+1}, \;\;\;\;\; \gamma(t) = -
 \frac{\gamma_0}{(t+1)^k}, \;\;\;\;\; \mu(t) = \mu_0 (t+1)^m ,
 \end{equation}
 where $d(t)$ is the lower bound of the diffusion coefficients, see
(\ref{c3}),
 $d_0, \gamma_0,$ and $ \mu_0$ are some positive constants, $k \geq
 1$ is a constant. Then inequality (\ref{c8}) takes the form:

 \begin{equation}
 \label{ap3}
  \alpha(t)  \leq  \mu_0^{q-1} (t+1)^{m(q-1)} \left(
 c(\Omega) \frac{d_0}{t+1} - \frac{\gamma_0}{(t+1)^k}
 - \frac{m}{t+1}  \right) , \;\; t \geq 0
 \end{equation}
 Let us assume that

 \begin{equation}
 \label{ap4}
 c(\Omega) d_0 >  \gamma_0 + m .
 \end{equation}
 If the above inequality holds, then the right-hand side of (\ref{ap3})
is
positive. This
 inequality gives a condition on the function $c_0(t)$, defined in
equation (\ref{c2}) and used in the
 definition of $\alpha(t)$ in (\ref{alpha}). If $m (q-1) < 1$, then
condition (\ref{ap3}) implies that
$c_0(t)$ should converge to zero as $t\to \infty$, if $m (q-1) > 1$, then
it $c_0(t)$ may grow, as $t$ grows, and still inequality  (\ref{ap3})
may be satisfied.

To satisfy assumption (\ref{c9}) one may choose
 \begin{equation}
 \label{in5}
 \mu_0=g(0)^{-1}.
 \end{equation}
If (\ref{ap3}), (\ref{ap4}), and (\ref{in5}) hold, then one may apply
 estimate (\ref{c10}) and obtain the following
 theorem.

 \medskip

 \noindent
 {\bf Theorem 3.2.}
 {\it If conditions (\ref{ap3a})-
(\ref{in5}) are satisfied,
 then the $L^2$ norm of the solution $u(x,t)$ of problem
 (\ref{m1})-(\ref{m3}) with the Dirichlet boundary condition
 admits the estimate

 $$ \|u(\cdot,t)\| \; \leq \; g(0) \; (t+1)^{-m}  \; , \;\;\; \forall t
 \geq 0 . $$
}


 \bigskip

 {\bf Boundedness of the solution.}

 Consider the case when  global asymptotic stability of the
 stationary solution may not hold. We wish to obtain an estimate of
 the solution of the evolution problem, which  yields
 stability in the sense of Lyapunov. We will illustrate the method  in
 the case of Neumann boundary condition.

 If the Neumann boundary condition holds, then, in contrast with the
Dirichlet boundary condition,  one has
 $c(\Omega)=0$ in equation (\ref{alpha}) and inequality (\ref{c5a}), so
one gets $\sigma(t) \equiv \gamma(t)$.

Let
 $$ \gamma(t) = \frac{\gamma_0}{1+t}, \;\;\;\;\; \mu(t) = \mu_0 (1+t)^m
, $$
 where $\gamma_0$ and $\mu_0$ are some positive constants.
 If $\gamma_0 > m$ and $c_0(t)$ is such that inequalities (\ref{c8})
and (\ref{c9}) hold, i.e.,
$$\mu_0=g(0)^{-1},$$
and
$$  Cc_0(t)\le \mu_0^{q-1}(1+t)^{m(q-1)}\frac {\gamma_0- m}{1+t},$$
then inequality (\ref{c10}) yields convergence at the rate
$O((1+t)^{-m})$.
 In this example $\gamma(t)$ is positive.

We can  consider the case when  $\gamma(t)$ is {\it negative}, but then
$\mu(t)$ has to be a decreasing function. For instance, assume that

 \begin{equation}
 \label{z1}
 \gamma(t) = -\frac{\gamma_0}{(t+1)^k}, \;\;\;\;\;
 \mu(t) = \mu_0  + \mu_1 (t+1)^{-\nu},
 \end{equation}
 where the constants $\gamma_0, \mu_0, \mu_1 > 0$ and  $\nu>0$.
 In this case, (\ref{c10}) yields boundedness of
 the solution for all $t\ge 0$, but the solution does not  converge
to zero.

 Inequality (\ref{c8}) takes the form:

 \begin{equation}
 \label{z2}
 \alpha(t) \leq \left( \mu_0  + \mu_1 (t+1)^{-\nu}
 \right)^{q-1} \left( \frac{\nu \mu_1 (1+t)^{-\nu-1} }{\mu_0  +
 \mu_1 (t+1)^{-\nu}} - \frac{\gamma_0 }{(t+1)^k} \right).
 \end{equation}
This inequality holds if, for example,
$$\nu+1\le k,$$
and, with $\alpha (t)= C c_0(t)$, the following inequality holds:
\begin{equation}
 \label{z2a}
 C(1+t)^{\nu +1} c_0(t) \leq  \mu_0^{q-1} \left( \frac{\nu\mu_1}
{\mu_0}   - \gamma_0 \right).
 \end{equation}
If (\ref{z2a}) holds, and
\begin{equation}
 \label{z22}
\mu(0)=\mu_0+\mu_1=[g(0)]^{-1},
\end{equation}
then inequality
(\ref{c10}) yields the following Theorem.

 \medskip

 \noindent
 {\bf Theorem 3.3.}
{\it  If conditions (\ref{z1})-(\ref{z2a}) hold, and $\nu+1\le k$,
  then the $L^2$ norm of the solution $u(x,t)$ of problem
 (\ref{m1})-(\ref{m3}) with the Neumann boundary condition
 satisfies the estimate

 $$ \|u(\cdot,t)\| \; \leq \;  [\mu (t)]^{-1}\le [\mu_0]^{-1}  \qquad
\forall t \geq 0 . $$

}
 \bigskip


 \subsection{Time-dependent Turing structures}

 Consider a reaction-diffusion system

 \begin{equation}
 \label{tur1}
 \frac{\partial u}{\partial t} = d_1(t) \frac{\partial^2 u}{\partial x^2} +
F(u,v,t) ,
 \end{equation}

 \begin{equation}
 \label{tur2}
 \frac{\partial v}{\partial t} = d_2(t) \frac{\partial^2 v}{\partial x^2} +
G(u,v,t)
 \end{equation}
 in the interval $0 < x < L$ with the boundary conditions

 \begin{equation}
 \label{tur3}
 u(0)=u(L)=0, \;\;\; v(0)=v(L)=0 .
 \end{equation}
 Reaction-diffusion systems describe various applied problems, for
example, biological problems.
 These systems are often considered in the case when the
coefficients
 and the nonlinearities do not depend explicitly on time. We
 introduce time dependence in order to describe variations
 of the environment (e.g., climate factors),  or to control the system
 behavior. For instance, if $u$ and $v$ are some
 concentrations, then the coefficients of mass diffusion and the
 reaction rates can depend on the temperature that can change in time due
to some external conditions, or the temperature can serve as a control
parameter.

 Suppose that $F(0,0,t) = G(0,0,t) = 0$ for all $t \geq 0$,
 that is $u=v=0$ is a stationary solution of problem (\ref{tur1})-(\ref{tur3}).
 This zero solution is also a stationary solution to the ODE system

\begin{equation}
 \label{tur4}
 \frac{d u}{d t} =  F(u,v,t) ,
 \end{equation}

 \begin{equation}
 \label{tur5}
 \frac{d v}{d t} =  G(u,v,t) .
 \end{equation}
 To simplify  calculations, let us assume that

 $$ F(u,v,t) = \phi(t) \; F^0(u,v), \;\;\; G(u,v,t) = \phi(t)
 \; G^0(u,v), \;\;\; d_i(t) = \phi(t) \; d_i^0, \;\; i=1,2 . $$

 Consider first the case where $\phi(t) \equiv 1$. Let us choose
 parameters in such a way that $u=v=0$ is  a stable
 solution of system (\ref{tur4})- (\ref{tur5}) but it is unstable
 as a solution of problem (\ref{tur1})-(\ref{tur3}).
 In this case, another solution, which is not homogeneous in
 space, can appear. This is so-called {\em Turing
 structure}, that is often discussed in relation with numerous biological
 applications (see, e.g., \cite{turing}, \cite{meinhardt}-\cite{othmer}).
 The Turing structure  provides one of the possible mechanisms of pattern
formation in biology.

 We assume that the solution $u=v=0$ of system (\ref{tur4}),
 (\ref{tur5}) is stable, and that the eigenvalues of the matrix

 $$ M = \left( \begin{array}{ccc}
 a &  & b \\
 c &  & d \end{array} \right) , $$
 have negative real parts. Here

 $$ a = F_u^0(0,0), \;\;\; b=F_v^0(0,0) , \;\;\; c=G_u^0(0,0),
\;\;\; d=G_v^0(0,0). $$
 In order to study stability of the zero  solution as a stationary
 solution of problem (\ref{tur1})-(\ref{tur3}), consider the
 linearized system

 \begin{equation}
 \label{lin1}
 \frac{\partial u}{\partial t} = d_1^0 \frac{\partial^2 u}{\partial x^2} +
a u + b v ,
 \end{equation}

 \begin{equation}
 \label{lin2}
 \frac{\partial v}{\partial t} = d_2^0 \frac{\partial^2 v}{\partial x^2} +
 c u + d v
 \end{equation}
 with the boundary conditions (\ref{tur3}).

 If one looks for the solutions of this system in the form

 $$ u = p_1 \sin (kx) e^{\lambda t} , \;\;\; v = p_2 \sin (kx) e^{\lambda t}  , $$
 then one obtains the following eigenvalue problem:

 $$ M(k) p = \lambda p , $$
 where
 $$ M(k) = \left( \begin{array}{ccc}
 a -d_1^0 k^2 &  & b \\
 c &  & d - d_2^0 k^2 \end{array} \right) , $$
 $p=(p_1,p_2)$, $M(0)=M$.
 Denote its eigenvalues by $\lambda_i(k), i=1,2$.
 The assumption Re $\lambda_i(0) < 0$, $i=1,2$ implies:

 \begin{equation}
 \label{tur11}
 a+d < 0, \;\;\;\; ad - bc > 0 .
 \end{equation}
 Furthermore,

 $$ \det M(k) = ad-bc - (a d_2^0 + d d_1^0) k^2 +d_1^0 d_2^0 k^4 , \;\;\;
 {\rm Tr} \; M(k) = a + d - (d_1^0 + d_2^0) k^2 . $$
 If $\det M(k) = 0$, then one eigenvalue of this matrix is
 negative and another one equals zero. Hence,
  system (\ref{lin1}), (\ref{lin2}) linearized about the solution
  $u=v=0$ has a zero eigenvalue. If, under a change of parameter,
  this eigenvalue crosses the origin, then a spatially
  inhomogeneous solution can bifurcate from it.

  Thus, equality $\det M(k) = 0$ determines the stability
  boundary of the solution $u=v=0$ and the condition of
  bifurcation of a spatially inhomogeneous solution.

 Let us verify that equality $\det M(k) = 0$  is
 compatible with inequalities (\ref{tur11}).
 If $a<0$ and $d<0$, then $\det M(k) > 0$.
 In order to have $\det M(k)=0$, let us assume that one of the
coefficients $a$ or $d$ is positive, but the sum $a+d$ is negative.
Let us assume, for instance, that $a>0$. The constants $b,c,d$ can be
chosen in
 such a way that inequalities (\ref{tur11}) are satisfied. For some
$a,b,c,d, k, d_2^0$
 fixed, we can increase $d_1^0$ so that the determinant
 of the matrix $M(k)$ becomes zero.

 \medskip

 Thus, if $\phi(t) \equiv 1$, then $u=v=0$ can be a stable solution
 of system (\ref{tur4}), (\ref{tur5}), but unstable as a solution
 of problem (\ref{tur1})-(\ref{tur3}). In this case, a stationary
 spatial structure can emerge and the solution of the evolution
 problem can converge to it.

 If $\phi(t) \not\equiv 1$, then the previous considerations do
 not allow us to conclude whether the solution of problem
 (\ref{tur1})-(\ref{tur3}) with a given initial condition
 converges to a trivial solution or to a spatially inhomogeneous
 solution.

 Let us use
 the method developed in Section 2 in order to study the behavior
 of solution of the time dependent reaction-diffusion system.
 We have

 $$ \sigma(t) = \phi(t) (d_0 c(\Omega) - \gamma_0) , $$
 where

 $$ d_0 = \min (d_1^0, d_2^0) , $$

 $$ a u_1^2 + (b+c) u_1 u_2 + d u_2^2 \leq \left(a + \frac12 (b+c) \right)
 u_1^2 + \left(d + \frac12 (b+c) \right) u_2^2 \leq \gamma_0 (u_1^2 + u_2^2),
 $$

 $$ \gamma_0 = \max \left( a + \frac12 (b+c), d + \frac12 (b+c)
 \right) . $$
 We obtain the following result.

 \medskip

 \noindent
 {\bf Theorem 3.4.}
 1. Assume that $d_0 c(\Omega) > \gamma_0$,
 $$ \phi(t) = \frac{\phi_0}{t+1}\; , \;\;\;\; \mu(t) = \mu_0
 (t+1)^m ,  \;\;\;\; \mu_0^{-1}=g(0). $$
 If $\phi_0 (d_0 c(\Omega) - \gamma_0) > m$,
 and $c_0(t)$ (see (\ref{c2})) and $\alpha(t)$ (see (\ref{alpha}))  are such that condition (\ref{c8}) is
satisfied, then

 $$ \| u(\cdot,t)\| \leq g(0) (t+1)^{-m} , \;\;\; t \geq 0 . $$
 2. Let $d_0 c(\Omega) < \gamma_0$,

 $$ \mu(t) = \mu_0 + \frac{\mu_1}{(t+1)^k} \; , \;\;\;\;
 (\mu_0 + \mu_1)^{-1} = g(0) , $$
 where $\mu_0$, $\mu_1$ and $k$ are some positive constants.

 If $c_0(t)$ and $\alpha(t)$ satisfy condition (\ref{c8}), then

 $$  \| u(\cdot,t)\| \leq \mu_0^{-1} \quad \forall t\ge 0.$$

 \bigskip

 The conclusion of this theorem follows from Lemma 2.1. The first part
 of the theorem gives a sufficient condition of convergence to the
 trivial solution. If this condition is not satisfied, then the
 solution can possibly converge to a spatially inhomogeneous solution.
 In this case,
 the second part of the theorem gives an estimate of the solution.



 \section{An estimate of the solution}

  \medskip

  \noindent
  {\bf Lemma 4.1.}
 Suppose that for some positive constant $M_1$ the following
 estimate holds:

 \begin{equation}
 \label{est1}
 |F(u,x,t)| \leq M_1, \;\;\; \forall u \in \mathbb R^n, \;\; x \in
 \Omega, \;\;  t \geq 0 .
 \end{equation}
 Then solution of problem (\ref{m1})-(\ref{m3}) with the Dirichlet
 boundary condition satisfies the estimate

 \begin{equation}
 \label{est2}
  \|u\|_{H^2(\Omega)} \leq M_2, \;\;\; t \geq 0 .
 \end{equation}

 \medskip

 \noindent
 {\bf Proof.}
 Each component $u_i$ of the solution satisfies the problem

 \begin{equation}
 \label{est3}
 \frac{\partial u_i}{\partial t} = d_i \Delta u_i + F_i(u,x,t) ,
 \end{equation}

 \begin{equation}
 \label{est4}
 u_i |_{\partial \Omega} = 0, \;\;\; u_i(x,0) = u_i^0(x) .
 \end{equation}
 We first obtain an estimate of the solution in the uniform
 norm. Let

 $$ v(x) = -a |x|^2 + b , $$
 where
 $$ - 2 n a d_i + M_1 \leq 0.  $$
 The ball with the radius $R=\sqrt{b/a}$ contains the domain
 $\Omega$, and
 $$ u_i^0(x) \leq v(x) , \;\;\; x \in \Omega . $$
 Such constants $a$ and $b$ can be chosen for any $M_1$ and any
 initial condition.

Then $v(x)$ is an upper solution of equation (\ref{est3}) and

 $$ u_i(x,t) \leq v(x), \;\;\; x \in \Omega, \;\; t \geq 0 . $$
The functions $u_i(x,t)$ are bounded from below as well. Thus,
estimate (\ref{est2}) follows from the known estimate (see, e.g.,
\cite{f})

 $$ \|u\|_{H^2(\Omega)} \leq K \left( \|F\|_{L^2(\Omega)} +
 \|u\|_{L^2(\Omega)} \right) ,$$
where $K>0$ is a constant independent on $u$.

 \hspace*{14cm}  $\Box$

 \medskip

  \noindent
  {\bf Lemma 4.2.}
 Suppose that

 \begin{equation}
 \label{est5}
 F_i(u,x,t) \leq 0, \;\;\; \forall u_i \geq u_i^*, \;\; x \in
 \Omega, \;\;  t \geq 0 ,
 \end{equation}
 for some constants $u_i^*$, $1\le i \le n$.
 Then solution of problem (\ref{m1})-(\ref{m3}) with the Neumann
 boundary condition satisfies estimate (\ref{est2}).

 \medskip

 \noindent
 {\bf Proof.}
 It is sufficient to note that any constant greater than $u_i^*$
 is an upper solution of equation (\ref{est3}).

 \hspace*{14cm}  $\Box$


\end{document}